\documentclass{amsart}
\usepackage{amsfonts,amsmath,amsthm,amssymb}
\newtheorem{theorem}{Theorem}

\newtheorem{proposition}{Proposition}
\begin{document}

\title{TORSION UNITS IN INTEGRAL GROUP RINGS\\ OF CONWAY SIMPLE GROUPS\\
\vspace{10pt}\Small{CIRCA preprint 2010/12}}

\dedicatory{The paper is dedicated to Professor Said Sidki on his 70th birthday.}

\author{V.A.~Bovdi}

\address{Institute of Mathematics, University of Debrecen,
P.O.  Box 12, H-4010 Debrecen, Hungary}

\author{A.B.~Konovalov}

\address{School of Computer Science, University of St Andrews,
North Haugh, St Andrews, Fife, KY16 9SX, Scotland}

\author{S.~Linton}

\address{School of Computer Science, University of St Andrews,
North Haugh, St Andrews, Fife, KY16 9SX, Scotland}

\begin{abstract}
Using the Luthar--Passi method, we investigate the possible
orders and partial augmentations
of torsion units of the normalized unit group of
integral group rings of Conway simple groups 
$\texttt{Co}_1$, $\texttt{Co}_2$ and $\texttt{Co}_3$.
\end{abstract}

\keywords{Zassenhaus conjecture; torsion unit; partial augmentation; 
integral group ring; Conway simple groups.}

\maketitle


Let $U(\mathbb Z G)$ be the unit group of the integral
group ring $\mathbb Z G$ of a finite group $G$, and
$V(\mathbb Z G)$ be its normalized unit group
$$
V(\mathbb ZG) = \Big \{ \; \sum_{g\in G}\alpha_g g \in U(\mathbb
ZG) \mid \sum_{g\in G}\alpha_g=1 \; \Big \}.
$$
The structure of $U(\mathbb Z G)$ is completely determined
by its normalized unit group since
$U(\mathbb Z G) = U(\mathbb Z) \times V(\mathbb Z G)$.
Throughout the paper (unless stated otherwise) any unit of $\mathbb Z G$
is always {\it normalized} and not equal to the identity element of $G$.

The following longstanding conjecture is due to H.~Zassenhaus 
(see \cite{Zassenhaus}):
\begin{itemize}
\item[]{\bf (ZC)} \qquad
every torsion unit $u \in V(\mathbb ZG)$ is conjugate within
the rational group algebra $\mathbb Q G$ to an element in $G$.
\end{itemize}
For finite simple groups the main tool for the investigation of the
Zassenhaus conjecture is the Luthar--Passi method, introduced in
\cite{Luthar-Passi} for the alternating group $A_{5}$ (for its
further applications, see also \cite{Bachle-Kimmerle} and \cite{Hertweck2}).

The conjecture {\bf (ZC)} is still open for all sporadic simple groups,
and for several of them results are available that either prove {\bf (ZC)} for some
orders or restrict possible partial augmentations
of torsion units. For some recent results on Mathieu, Janko,
Higman-Sims, McLaughlin, Held, Rudvalis, Suzuki and O'Nan simple groups
 we refer to 
\cite{Bovdi-Grishkov-Konovalov-HeON,Bovdi-Jespers-Konovalov-Janko,
Bovdi-Konovalov-M11, Bovdi-Konovalov-McL, 
Bovdi-Konovalov-M23, 
Bovdi-Konovalov-Ru, Bovdi-Konovalov-HS,
Bovdi-Konovalov-Linton-M22, Bovdi-Konovalov-Marcos-Suz, Bovdi-Konovalov-Siciliano-M12}. Here we continue these investigations for 
the Conway simple groups $\texttt{Co}_1$, $\texttt{Co}_2$ and $\texttt{Co}_3$.

Let $G$ be a finite group.
Denote by $\mathcal{C} =\{ C_{1}, C_{k_1 t_1}, \ldots, C_{k_s t_r} \}$
the collection of all conjugacy classes of $G$, where $C_{1}=\{ 1\}$,
and $C_{k t}$ denote the conjugacy class with representatives
of order $k$, labelled by the distinguishing letter $t$
(throughout the paper we use the ordering of
conjugacy classes as used in the GAP Character Table Library). Suppose
$u=\sum \alpha_g g \in V(\mathbb Z G)$ is a non-trivial torsion unit.
The partial augmentation of $u$ with respect to $C_{nt}$ is defined
as $\nu_{nt}=\nu_{nt}(u)=\sum_{g\in C_{nt}} \alpha_{g}$.

The criterion for {\bf ZC} can be formulated in terms of the vanishing of 
partial augmentations of torsion units (see Proposition \ref{P:MRSW} below). 
Therefore, it is useful to know for each possible
order of a torsion unit in $V(\mathbb Z G)$, which combinations of partial 
augmentations may arise. Such an answer is provided by our main results.

\begin{theorem}\label{T:1}
Let $G$ denote the Conway simple group $\texttt{Co}_3$. Let $u$
be a torsion unit of $V(\mathbb ZG)$ of order $|u|$ with the tuple 
$\nu$ of length $42$ containing partial augmentations 
for all conjugacy classes of $G$. The following properties hold.

\begin{itemize}

\item[{\rm (i)}] There are no elements of order
$33$, $46$, $55$, $69$, $77$, $115$, $161$, $253$ in $V(\mathbb ZG)$.
Equivalently, if $|u| \not \in \{ 28$, $35$, $36$, $40$, $42$, $44$,
$45$, $56$, $60$, $63$, $70$, $72$, $84$, $88$, $90$, $105$, $120$,
$126$, $140$, $168$, $180$, $210$, $252$, $280$, $315$, $360$, $420$,
$504$, $630$, $840$, $1260$, $2520 \}$,
then $|u|$ is the order of some element $g \in G$.

\item[{\rm (ii)}] If $|u|=7$, then $u$ is rationally conjugate to some $g\in G$.

\item[{\rm (iii)}] If $|u|=2$, then \quad $\nu_{kx}=0$ \; for \; $kx \not\in \{ 2a, 2b \}$ and
$$(\nu_{2a}, \nu_{2b}) \in \{ \; (-2,3), \; (-1,2), \; (0,1), \; (1,0), \; (2,-1), \; (3,-2) \; \}.$$

\item[{\rm (iv)}] If $|u|=3$, all partial augmentations of $u$ 
are zero except possibly $\nu_{3a}$,$\nu_{3b}$,$\nu_{3c}$ and the triple
$(\nu_{3a},\nu_{3b},\nu_{3c})$ is one of those given in Appendix A.

\item[{\rm (v)}] If $|u|=5$, then \quad $\nu_{kx}=0$ \; for \; $kx \not\in \{ 5a, 5b \}$ and
$$(\nu_{5a}, \nu_{5b}) \in \{ \; (-4,5), \; (-3,4), \; (-2,3), \; (-1,2), \; (0,1), \; (1,0) \; \}.$$

\item[{\rm (vi)}] If $|u|=11$, then all partial augmentations of $u$ are zero
except possibly $\nu_{11a}$,$\nu_{11b}$, and the pair $(\nu_{11a},\nu_{11b})$
is one of
$$
\{ \; (\nu_{11a},\nu_{11b}) \; \mid \; -11 \le  \nu_{11a} \le 12, \quad \nu_{11a}+\nu_{11b}=1 \; \}.
$$

\item[{\rm (vii)}] If $|u|=23$, then all partial augmentations of $u$ are zero
except possibly $\nu_{23a}$,$\nu_{23b}$, and the pair $(\nu_{23a},\nu_{23b})$
is one of
$$
\{ \; (\nu_{23a},\nu_{23b}) \; \mid \; -5 \le  \nu_{23a} \le 6, \quad \nu_{23a}+\nu_{23b}=1 \; \}.
$$

\item[{\rm (viii)}] If $|u|=35$, then all partial augmentations of $u$ are zero
except possibly $\nu_{5a}$, $\nu_{5b}$, $\nu_{7a}$,
and the triple $( \nu_{5a}, \nu_{5b}, \nu_{7a} )$
is one of \; $\{\, (3,12,-14), \, (4,11,-14) \, \}$.

\end{itemize}

\end{theorem}

\begin{theorem}\label{T:2}
Let $G$ denote the Conway simple group $\texttt{Co}_2$. Let  $u$
be a torsion unit of $V(\mathbb ZG)$ of order $|u|$ 
with the tuple $\nu$ of length $60$ containing partial augmentations 
for all conjugacy classes of $G$. The following properties hold.

\begin{itemize}

\item[{\rm (i)}] There are no elements of order
$21$, $22$, $33$, $46$, $55$, $69$, $77$, $115$, $161$, $253$ in $V(\mathbb ZG)$.
Equivalently, if $|u| \not \in \{ 36, 40, 45, 48, 56, 60, 70, 72, 80, 90, 112, 120$, 
$140, 144, 180, 240, 280, 360, 560, 720 \}$, then $|u|$ is the order of some $g \in G$.

\item[{\rm (ii)}] If $|u| \in \{ 7, 11 \}$, then $u$ is rationally conjugate to some $g\in G$.

\item[{\rm (iii)}] If $|u|=2$, all partial augmentations of $u$ are zero
except possibly $\nu_{2a}$,$\nu_{2b},\nu_{2c}$, and the triple $(\nu_{2a},\nu_{2b},\nu_{2c})$
is one of those given in Appendix B.

\item[{\rm (iv)}] If $|u|=3$, then \quad $\nu_{kx}=0$ \; for \; $kx \not\in \{ 3a, 3b \}$ and
$$(\nu_{3a}, \nu_{3b}) \in \{ \; (-2,3), \; (-1,2), \; (0,1), \; (1,0) \; \}.$$

\item[{\rm (v)}] If $|u|=5$, then \quad $\nu_{kx}=0$ \; for \; $kx \not\in \{ 5a, 5b \}$ and
$$(\nu_{5a}, \nu_{5b}) \in \{ \; (-4,5), \; (-3,4), \; (-2,3), \; (-1,2), \; (0,1), \; (1,0) \; \}.$$

\item[{\rm (vi)}] If $|u|=23$, then all partial augmentations of $u$ are zero
except possibly $\nu_{23a}$,$\nu_{23b}$, and the pair $(\nu_{23a},\nu_{23b})$
is one of
$$
\{ \; (\nu_{23a},\nu_{23b}) \; \mid \; -32 \le  \nu_{23a} \le 33, \quad \nu_{23a}+\nu_{23b}=1 \; \}.
$$

\item[{\rm (vii)}] If $|u|=35$, then all partial augmentations of $u$ are zero
except possibly $\nu_{5a}$, $\nu_{5b}$, $\nu_{7a}$,
and the triple $( \nu_{5a}, \nu_{5b}, \nu_{7a} )$
is one of \; $\{\, (3,12,-14), \, (4,11,-14) \, \}$.

\end{itemize}

\end{theorem}

\begin{theorem}\label{T:3}
Let $G$ denote the Conway simple group $\texttt{Co}_1$. Let  $u$
be a torsion unit of $V(\mathbb ZG)$ of order $|u|$ 
with the tuple $\nu$ of length $101$ containing partial augmentations 
for all conjugacy classes of $G$. The following properties hold.

\begin{itemize}

\item[{\rm (i)}] There are no elements of order
$46$, $69$, $77$, $91$, $115$, $143$, $161$, $253$ and $299$ in $V(\mathbb ZG)$.
Equivalently, if $|u| \not \in \{ 55, 65, 110, 130, 165, 195, 220, 260, 330, 390$,
$440, 495, 520, 585, 660, 780, 880, 990, 1040, 1170, 1320, 1560, 1980, 2340, 2640$, 
$3120, 3960, 4680, 7920, 9360 \}$, then $|u|$ is the order of some $g \in G$.

\item[{\rm (ii)}] If $|u| \in \{ 11, 13 \}$, then $u$ is rationally conjugate to some $g\in G$.

\item[{\rm (iii)}] If $|u|=7$, then all partial augmentations of $u$ are zero
except possibly $\nu_{7a}$ and $\nu_{7b}$, and the pair $(\nu_{7a},\nu_{7b})$
is one of $$
\{ \; (\nu_{7a},\nu_{7b}) \quad \mid \quad \nu_{7a} + \nu_{7b} =1, \quad -7 \le \nu_{7a} \le 39 \; \}.$$

\item[{\rm (iv)}] If $|u|=23$, then all partial augmentations of $u$ are zero
except possibly $\nu_{23a}$ and $\nu_{23b}$, and the pair $(\nu_{23a},\nu_{23b})$
is one of 
$$
\{ \; (\nu_{23a},\nu_{23b}) \quad \mid \quad \nu_{7a} + \nu_{7b} = 1, \quad -29293 \le \nu_{23a} \le 29294 \; \}.
$$

\item[{\rm (v)}] If $|u|=55$, then all partial augmentations of $u$ are zero
except possibly $\nu_{5a}$, $\nu_{5b}$, $\nu_{5c}$ and $\nu_{11a}$, 
and the tuple $( \nu_{5a}, \nu_{5b}, \nu_{5c}, \nu_{11a} )$ is one of
\[
\begin{aligned}
( -2, 2, -10, 11 ), \;  ( -2, 3, -11, 11 ), \;  ( -1, -2, -7, 11 ), \;  ( -1, -1, -8, 11 &),  \\
( -1, 0, -9, 11 ), \;  ( -1, 1, -10, 11 ), \;  ( 0, -6, -4, 11 ), \;  ( 0, -5, -5, 11 &),  \\
( 0, -4, -6, 11 ), \;  ( 0, -3, -7, 11 ), \;  ( 0, -2, -8, 11 ), \;  ( 0, -1, -9, 11 &), \\
( 0, 0, -10, 11 ), \;  ( 1, -8, -3, 11 ), \;  ( 1, -7, -4, 11 ), \;  ( 1, -6, -5, 11 &),  \\
( 1, -5, -6, 11 ), \;  ( 1, -4, -7, 11 ), \;  ( 1, -3, -8, 11 ), \;  ( 1, -2, -9, 11 &),  \\
( 2, -9, -3, 11 ), \;  ( 2, -8, -4, 11 ), \;  ( 2, -7, -5, 11 ), \;  ( 2, -6, -6, 11 &),  \\
( 2, -5, -7, 11 ), \;  ( 2, -4, -8, 11 ), \;  ( 2, -3, -9, 11 ), \;  ( 3, -11, -2, 11 &),  \\
( 3, -10, -3, 11 ), \;  ( 3, -9, -4, 11 ), \;  ( 3, -8, -5, 11 ), \;  ( 3, -7, -6, 11 &),  \\ 
( 3, -6, -7, 11 ), \;  ( 4, -12, -2, 11 ), \;  ( 4, -11, -3, 11 ), \;  ( 4, -10, -4, 11 &).
\end{aligned}
\]

\item[{\rm (vi)}] If $|u|=65$, then all partial augmentations of $u$ are zero
except possibly $\nu_{5a}$, $\nu_{5b}$, $\nu_{5c}$ and $\nu_{13a}$, 
and the tuple $( \nu_{5a}, \nu_{5b}, \nu_{5c}, \nu_{13a} )$ is one of
\[
\begin{aligned}
( -3, 2, -24, 26 ), \;  ( -2, -2, -21, 26 ), \;  ( -2, -1, -22, 26 ), \;  ( -2, 0, -23, 26 & ), \\
( -1, -3, -21, 26 ), \;  ( -1, -2, -22, 26 ), \;  ( -1, -1, -23, 26 ), \;  ( 5, -4, 39, -39 & ), \\
( 5, -3, 38, -39 ), \;  ( 6, -7, 41, -39 ), \;  ( 6, -6, 40, -39 ), \;  ( 6, -5, 39, -39 & ), \; \\
( 7, -8, 41, -39 ), \;  ( 7, -7, 40, -39 & ).
\end{aligned}
\]

\end{itemize}

\end{theorem}


For the determination of possible orders of torsion units in $V(\mathbb ZG)$
first of all we start with the following well-known bound.

\begin{proposition}[\cite{Cohn-Livingstone}]\label{P:order(Cohn-Livingstone)}
The order of a torsion element $u\in V(\mathbb ZG)$ divides $exp(G)$.
\end{proposition}

Moreover, the partial augmentations of torsion units are also bounded.

\begin{proposition}[see \cite{Hales-Luthar-Passi}]\label{P:Hales-Luthar-Passi}
Let $C_1,\dots,C_n$ be conjugacy classes of a finite group $G$.
Let $u$ be a torsion unit in $V(\mathbb ZG)$ and
$\nu_i(u)$ denote the partial augmentation of $u$ with respect to the
conjugacy class $C_i$. Then \; $\nu_{i}(u)^2 \le |C_i|$ \; and, moreover,
$$
\sum_{i=1}^{n} \frac{\nu_{i}(u)^2}{|C_i|} \le 1.
$$
\end{proposition} 

The following result allows a reformulation of
the Zassenhaus conjecture in terms of the vanishing of
partial augmentations of torsion units.

\begin{proposition}[see \cite{Luthar-Passi}]\label{P:MRSW}
Let $u\in V(\mathbb Z G)$
be of order $k$. Then $u$ is conjugate in $\mathbb
QG$ to an element $g \in G$ if and only if for
each $d$ dividing $k$ there is precisely one
conjugacy class $C$ with partial augmentation
$\varepsilon_{C}(u^d) \neq 0 $.
\end{proposition}

The next result is a reformulation of the Proposition 3.1 \cite{Hertweck2}
(which was originally proved for a group ring over an arbitrary Dedekind 
ring of characteristic zero) for the case of integral group rings. 
It restricts possible values of some partial augmentations of torsion units.

\begin{proposition}[see \cite{Hertweck2}, Proposition 3.1]\label{P:Hertweck}
Let $G$ be a finite
group and let $u$ be a torsion unit in $V(\mathbb ZG)$ of order $k$.
If $x$ is an element of $G$ whose order does not divide $k$,
then $\varepsilon_x(u)=0$.
\end{proposition}

The basis of the Luthar--Passi method which produces
further restrictions on possible orders of torsion 
units and their partial augmentations is the following.

\begin{proposition}[see \cite{Hertweck2,Luthar-Passi}]\label{P:Luthar-Passi}
 Let either $p=0$ or $p$ a prime divisor of $|G|$. Suppose
that $u\in V( \mathbb Z G) $ has finite order $k$ and assume $k$ and
$p$ are coprime in case $p\neq 0$. If $z$ is a complex primitive $k$-th root
of unity and $\chi$ is either a ordinary character or a $p$-Brauer
character of $G$, then for every integer $l$ we define 
\begin{equation}\label{E:2}
\mu_l(u,\chi, p ) =
\textstyle\frac{1}{k} \sum_{d|k}Tr_{ \mathbb Q (z^d)/ \mathbb Q }
\bigl( \chi(u^d)z^{-dl} \bigr).
\end{equation}
Then $\mu_l(u,\chi, p )$ is a non-negative integer not greater than $\text{deg}(\chi)$.
\end{proposition}

Note that if $p=0$, we will use the notation $\mu_l(u,\chi,*)$ for
$\mu_l(u,\chi,  0)$.


Let $u$ be a normalized unit of order $k$, where $k$
divides $exp(G)$ by Proposition 
\ref{P:order(Cohn-Livingstone)}.
From the Berman--Higman Theorem (see \cite{Artamonov-Bovdi})
one knows that $\text{tr}(u) = \nu_1=0$, so 
\begin{equation}\label{E:1}
\sum_{C_{nt}\in \mathcal{C}\setminus C_1} \nu_{nt}=1.
\end{equation}
On the next step we apply Proposition \ref{P:Hertweck} 
for every appropriate prime $p$, such that $k=p^m t$, where $(p,t)=1$,
to eliminate partial augmentations
of conjugacy classes of elements of $G$ with representatives
of order $p^n s$, where $(p,s)=1$ and $n>m$. If after this step for 
torsion units of some order $k$ we have only one non-zero partial augmentation, 
then {\bf ZC} holds for this order by Proposition \ref{P:MRSW}. 

Otherwise, we have to produce and solve a system of constraints.
For the unit $u$ of order $k$ we denote by
$\nu_1^{(k)}, \dots, \nu_n^{(k)}$ its non-vanishing partial augmentations 
for conjugacy classes $C_{q_1}, \dots, C_{q_n}$ (we will also omit the upper 
index and denote them by $\nu_1, \dots, \nu_n$ for the clarity of notation).
Let $d_1, \dots, d_s$ be the set of all non-negative integers dividing $k$,
where $d_i > 1$ and $d_s=k$. Furthermore, let $k_i = k/d_i$
and let $\nu_1^{(k_i)}, \dots, \nu_{n_i}^{(k_i)}$ be the non-vanishing partial 
augmentations for elements of order $k_i$, corresponding to conjugacy classes
$C_{q_1}^{(k_i)}, \dots, C_{q_{n_i}}^{(k_i)}$. Then the right-hand side in
(\ref{E:2}) from Proposition \ref{P:Luthar-Passi} formula may be written as
\[
\begin{aligned}
\textstyle\frac{1}{k} \sum_{d|k} Tr_{ \mathbb Q (z^d)/ \mathbb Q } 
\bigl( \chi(u^d)z^{-dl} \bigr) = \textstyle\frac{1}{k} \bigl( \,
& Tr_{\mathbb Q(z)/\mathbb Q}\bigl(\chi(u)z^{-l}\bigr) \, + \\
& Tr_{\mathbb Q(z^{d_1})/\mathbb Q}\bigl(\chi(u^{d_1})z^{-d_{1}l}\bigr) \, + \, \cdots \, + \\
& Tr_{\mathbb Q(z^{d_i})/\mathbb Q}\bigl(\chi(u^{d_i})z^{-d_{i}l}\bigr) \, + \, \cdots \, + \\
& Tr_{\mathbb Q(z^{d_{s-1}})/\mathbb Q}\bigl(\chi(u^{d_{s-1}})z^{-d_{s-1}l}\bigr) \, + \, \chi(1) \, \bigr),
\end{aligned}
\]
where the summand $\chi(1)$ comes from $d_s=k$.

Clearly, $\chi(u)=\sum_{j=1}^n \; \chi(h_j) \; \nu_j$ and 
$\chi(u^{d_i})=\sum_{j=1}^{n_i} \; \chi(h_j) \; \nu_j^{(k_i)}$
for any character $\chi$,
where $h_j$ is a representative  of the conjugacy class $C_{j}$,
and $\nu_j^{(k_i)}$ is the partial augmentation for the conjugacy 
class $C_j$ for an element $u^{d_i}$ of order $k_i = d/d_i$.

Since the trace is a linear mapping, this gives us $\mu_l(u,\chi,p)$ 
as a linear combination of corresponding partial augmentations:
\[
\begin{aligned}
\mu_l(u,\chi,p) = \textstyle\frac{1}{k} \big( & c_1 \nu_1 + \dots + c_n \nu_n \; + \\
& c_1^{(k_1)} \nu_1^{(k_1)} + \dots + c_n^{(k_1)} \nu_n^{(k_1)} \; + \; \cdots \; + \\
& c_1^{(k_i)} \nu_1^{(k_i)} + \dots + c_n^{(k_i)} \nu_n^{(k_i)} \; + \; \cdots \; + \\
& c_1^{(k_{s-1})} \nu_1^{(k_{s-1})} + \dots + c_n^{(k_{s-1})} \nu_n^{(k_{s-1})} \; + \; \chi(1) \big) \ge 0.
\end{aligned}
\]
Since  all the trace values must lie in $\mathbb{Q}$,  we may be able to deduce at this stage that some more partial augmentations must be zero,
when the corresponding character values are irrational.

Now to form the constraint satisfaction problem (CSP) for units of order $k$ we put together:
all inequalities for $\mu_l(u,\chi,p_i)$ for units of order $k$
for all possible $0 \le l < k$, characters $\chi$ and $p_i$;
similarly produced on earlier steps systems of inequalities with 
indeterminates $\nu_1^{(k_i)}, \dots, \nu_{n}^{(k_i)}$ for units of order $k_i$;
equation $\nu_1 + \dots + \nu_n = 1$
and equations $\nu_1^{(k_i)} + \dots + \nu_{n_i}^{(k_i)} =1$ for every order $k_i$. 
Now, if this CSP has no solutions, this can be seen immediately, and this approach
is much more efficient than enumerating all cases determined by possible partial 
augmentations for units of orders $k_i$, used, for example, in 
\cite{Bovdi-Konovalov-M11, Bovdi-Konovalov-Siciliano-M12}.


Proposition \ref{P:Luthar-Passi} may be reformulated for elements of order
$st$. Let $s$ and $t$ be two primes such that $G$ contains no
element of order $st$, and let $u$ be a normalized torsion unit of
order $st$. We denote by $\nu_k$ the sum of partial augmentations
of $u$ with respect all conjugacy classes of elements of order $k$
in $G$, i.e. $\nu_2 = \nu_{2a}+\nu_{2b}$, etc. Then by (\ref{E:1})
and Proposition \ref{P:Hertweck} we obtain that $\nu_s + \nu_t = 1$ and
$\nu_k = 0$ for $k \notin \{ s, t \}$. For each character $\chi$
of $G$ (an ordinary character or a Brauer character in
characteristic not dividing $st$) that is constant on the
elements of order $s$ and constant on the elements of order $t$,
we have $\chi(u) = \nu_s \chi(C_s) + \nu_t \chi(C_t)$, where
$\chi(C_r)$ denote the value of the character $\chi$ on any
element of order $r$ of $G$.

Let $s$ and $t$ be two primes dividing $|G|$,
and let $\chi$ be an ordinary or $p$-Brauer character of $G$
for $p$ not dividing $st$. Then $\chi$ is called 
a $(s,t)${\it -constant character}, if $\chi$ is constant on all elements of 
order $s$ and constant on all elements of order $t$.
\vspace{5pt}

From Proposition \ref{P:Luthar-Passi} we obtain that the values
\begin{equation}\label{E:3}
\begin{aligned}
\mu_l(u, \chi) =
\textstyle \frac{1}{st} \Bigl( \; \chi(1)
& + Tr_{\mathbb Q(z^s)/ \mathbb Q} \bigl( \chi(u^s) z^{-sl} \bigr) \\
& + Tr_{\mathbb Q(z^t)/ \mathbb Q} \bigl( \chi(u^t) z^{-tl} \bigr)
  + Tr_{\mathbb Q(z)/ \mathbb Q} ( \chi(u) z^{-l} \bigr) \; \Bigr)
\end{aligned}
\end{equation}
are nonnegative integers, and if $\chi$ is $(s,t)$-constant character then we get
\begin{equation}\label{E:4}
  \mu_l(u, \chi) =
     \textstyle \frac{1}{st} \left( m_1 +  \nu_s m_s + \nu_t m_t \right),
\end{equation}
where
\begin{equation}\label{E:5}
\begin{aligned}
m_1 & = \chi(1) + \chi(C_t) \, Tr_{\mathbb Q(z^s)/\mathbb Q}( z^{-sl} )
              + \chi(C_s) \, Tr_{\mathbb Q(z^t)/\mathbb Q}( z^{-tl} ), \\
m_s & = \chi(C_s) \, Tr_{\mathbb Q(z)/\mathbb Q}( z^{-l} ), \qquad
m_t   = \chi(C_t) \, Tr_{\mathbb Q(z)/\mathbb Q}( z^{-l} ).
\end{aligned}
\end{equation}
Since Proposition \ref{P:Luthar-Passi} and its reformulation are valid for any 
character (not necessarily irreducible), we are interested in a systematic 
search for $(s,t)$-constant characters
that are capable of producing new constraints 
on partial augmentations. For example, if we have only two conjugacy 
classes of elements of order $k$, namely $C_{ka}$ and $C_{kb}$, then 
if there are two characters $\chi_{1}$ and $\chi_{2}$ such that 
$$\chi_{1}(ka)-\chi_{1}(kb) = \chi_{2}(kb)-\chi_{2}(ka),$$ then for the 
character $\chi=\chi_1+\chi_2$ we will have that $\chi(ka)=\chi(kb)$.

If we have two $(s,t)$-constant characters $\chi_{1}$ and $\chi_{2}$,
then $\chi_{1}+\chi_{2}$ can not give us any further 
restrictions on partial augmentations, as it is shown by the following.

\begin{proposition}
Let either $p=0$ or $p$ a prime divisor of $|G|$. Suppose
that $u\in V( \mathbb Z G) $ has finite order $k$ and assume $k$ and
$p$ are coprime in case $p\neq 0$. If $z$ is a complex primitive $k$-th root
of unity and $\chi_1$, $\chi_2$ are both either classical characters or 
$p$-Brauer characters of $G$, then $\mu_l(u,\chi_1+\chi_2, p )$ 
is a non-negative integer whenever both $\mu_l(u,\chi_1, p )$
and $\mu_l(u,\chi_2, p )$ are non-negative integers.
\end{proposition}

Indeed, put $\xi=\chi_1+\chi_2$. It is easy to check that
\[
\begin{aligned}
\mu_l(u,\xi, p ) & = 
\textstyle\frac{1}{k} \sum_{d|k}Tr_{ \mathbb Q (z^d)/ \mathbb Q }
\{(\chi_1+\chi_2)(u^d)z^{-dl}\}
\\
& =\textstyle\frac{1}{k} \sum_{d|k}Tr_{ \mathbb Q (z^d)/ \mathbb Q }
\{\chi_1(u^d)z^{-dl}\} +
\textstyle\frac{1}{k} \sum_{d|k}Tr_{ \mathbb Q (z^d)/ \mathbb Q }
\{\chi_2(u^d)z^{-dl}\}
\\
&=\mu_l(u,\chi_1, p ) + \mu_l(u,\chi_2, p ).
\end{aligned}
\]
Thus, the task is to find all $(s,t)$-constant characters that can
not be represented as a sum of other $(s,t)$-constant characters.
We will call such characters $(s,t)${\it -irreducible characters}.
The search can be performed by analyzing relative differences between
values of irreducible characters on all conjugacy classes of the 
given order (see example in the proof for units of order 
35 from Theorem \ref{T:1}).

\section{Proof of Theorem \ref{T:1} }

Let $G \cong \texttt{Co}_3$. It is
well known \cite{Conway,GAP} that 
$|G| = 2^{10} \cdot 3^{7} \cdot 5^{3} \cdot 7 \cdot 11 \cdot 23$  
and
$ exp(G) = 2^{3} \cdot 3^{2} \cdot 5 \cdot 7 \cdot 11 \cdot 23$.  
The ordinary
and  $p$-Brauer character tables of $G$ for
$p\in\{2,3,5,7,11,23\}$ can be found
using the computational algebra system GAP \cite{GAP}, which
derives its data from \cite{AFG,ABC}. For characters and
conjugacy classes we will use throughout the paper the same
notation, including indexation, as used in the GAP Character Table
Library.

The group $G$ possesses elements of orders
2,  3,  4,  5,  6,  7,  8,  9,  10,  11,  12,  14,  15,  18,  20,  21,  22,  23,  24 and 30.
We begin with units of prime orders: 2, 3, 5, 7, 11, 23.
We do not give here our results for the remaining cases of torsion units of orders
4, 6, 8, 9, 10, 12, 14, 15, 18, 20, 21, 22, 24 and 30 because they are rather complex.  For example, using our implementation of the Luthar--Passi method, which we
intend to make available in the GAP package LAGUNA \cite{LAGUNA}, together
with constraint solvers MINION \cite{MINION} and ECLiPSe \cite{ECLiPSe},
we can compute 510 possible
cases for partial augmentations $(\nu_{2a}, \nu_{2b}, \nu_{4a}, \nu_{4b})$ for torsion units
of order 4 and five possible cases for partial augmentations $(\nu_{2a}, \nu_{2b}, \nu_{7a}, \nu_{14a})$
for torsion units of order 14. 
To complete the proof, we will investigate units of orders which do not appear in $G$.

%
%

\noindent $\bullet$ Let $|u|=7$. Using
Proposition \ref{P:Hertweck} we obtain that all partial
augmentations except one are zero. Thus by Proposition \ref{P:MRSW}
part (ii) of  Theorem \ref{T:1} is proved.

%
%

\noindent $\bullet$ Let $|u|=2$.
By (\ref{E:1}) and Proposition \ref{P:Hertweck} we get
$\nu_{2a}+\nu_{2b}=1$. Put $t_1=7\nu_{2a}-\nu_{2b}$. 
Now using Proposition \ref{P:Luthar-Passi}
we obtain the following system of inequalities
$$
\mu_{0}(u,\chi_{2},*) = \textstyle \frac{1}{2} (t_1 + 23) \geq 0; \quad 
\mu_{1}(u,\chi_{2},*) = \textstyle \frac{1}{2} (-t_1 + 23) \geq 0, \\ 
$$
which has only six solutions listed in part (iii) of Theorem \ref{T:1}.

%
%

\noindent $\bullet$ Let $u$ be a unit of order $3$. By (\ref{E:1})
and Proposition \ref{P:Hertweck} we get $\nu_{3a}+\nu_{3b}+\nu_{3c}=1$.
Put $t_1=4 \nu_{3a} - 5 \nu_{3b} +  \nu_{3c}$,
$t_2=10 \nu_{3a} + 10 \nu_{3b} +  \nu_{3c}$ and
$t_3=14 \nu_{3a} + 5 \nu_{3b} + 2 \nu_{3c}$.
Using Proposition \ref{P:Luthar-Passi} we obtain
the following system of inequalities
\[
\begin{aligned}
\mu_{0}(u,\chi_{2},*) & = \textstyle \frac{1}{3} (-2t_1 + 23) \geq
0; \qquad
\mu_{1}(u,\chi_{2},*)  = \textstyle \frac{1}{3} ( t_1+ 23) \geq 0; \\ 
\mu_{0}(u,\chi_{3},*) & = \textstyle \frac{1}{3} (2t_2 + 253) \geq
0; \qquad
\mu_{1}(u,\chi_{3},*)  = \textstyle \frac{1}{3} (-t_2 + 253) \geq 0; \\ 
\mu_{0}(u,\chi_{3},2)  & = \textstyle \frac{1}{3} (2t_3+ 230) \geq 0; \qquad 
\mu_{1}(u,\chi_{3},2) = \textstyle \frac{1}{3} (-t_3+ 230) \geq 0,
\end{aligned}
\]
which has only 155 solutions listed in the Appendix A
such that all $\mu_{i}(u,\chi_{j})$ are non-negative integers.

%
%

\noindent $\bullet$ Let $u$ be a unit of order $5$. By (\ref{E:1})
and Proposition \ref{P:Hertweck} we get $\nu_{5a}+\nu_{5b}=1$. Put $t_1=2
\nu_{5a} - 3 \nu_{5b}$. Using Proposition \ref{P:Luthar-Passi} we obtain the
system of two inequalities
\[
\mu_{0}(u,\chi_{2},*) = \textstyle \frac{1}{5} (-4t_1 + 23) \geq
0; \qquad \mu_{1}(u,\chi_{2},*) = \textstyle \frac{1}{5} (t_1 +
23) \geq 0,
\]
which has only six solutions listed in part (v) of Theorem \ref{T:1}
such that all $\mu_{i}(u,\chi_{2},*)$ are non-negative integers.

%
%

\noindent $\bullet$ Let $u$ be a unit of order $11$. By
(\ref{E:1}) and Proposition \ref{P:Hertweck} we get
$\nu_{11a}+\nu_{11b}=1$. Put $t_1=6 \nu_{11a} - 5 \nu_{11b}$ and
$t_2=-5 \nu_{11a} + 6 \nu_{11b}$ (observe that $t_2=1-t_1$). Now 
$$
\mu_{1}(u,\chi_{3},3) = \textstyle \frac{1}{11} (t_1 + 126) \geq 0; \qquad
\mu_{2}(u,\chi_{3},3)  = \textstyle \frac{1}{11} (-t_1 + 127) \geq 0, \\ 
$$
and this system has only 24 solutions listed in part (vi) of Theorem \ref{T:1}
such that all $\mu_{i}(u,\chi_{j},*)$ are non-negative integers.

%
%

\noindent $\bullet$ Let $u$ be a unit of order $23$. By
(\ref{E:1}) and Proposition \ref{P:Hertweck} we get
$\nu_{23a}+\nu_{23b}=1$. Put $t_1=12 \nu_{23a} - 11 \nu_{23b}$. 
Now using Proposition \ref{P:Luthar-Passi} we obtain the system of  inequalities
$$
\mu_{1}(u,\chi_{3},3) = \textstyle \frac{1}{23} (t_1 + 126) \geq 0; \qquad\;\;
\mu_{5}(u,\chi_{3},3) = \textstyle \frac{1}{23} (-t_1 + 127) \geq 0, \\ 
$$
which has only 12 solutions listed in part (vii) of Theorem \ref{T:1}
such that all $\mu_{i}(u,\chi_{j},*)$ are non-negative integers.

Now we consider orders which do not appear in $G$.

%
%

\noindent $\bullet$ Let $|u|=33$.
For these units we consider partial augmentations 
$\nu_{3a}$, $\nu_{3b}$, $\nu_{3c}$, $\nu_{11a}$ and $\nu_{11b}$.
Since $|u^3|=11$ and $|u^{11}|=3$, by
Proposition \ref{P:Luthar-Passi} we obtain the system of inequalities with 5 more variables
$\nu_{11a}^{(3)}$, $\nu_{11b}^{(3)}$, $\nu_{3a}^{(11)}$, $\nu_{3b}^{(11)}$ and $\nu_{3c}^{(11)}$.
Replacing $\nu_{11a}^{(3)}$ and $\nu_{11b}^{(3)}$ by their numerical values from
part (vi) of Theorem \ref{T:1}, we get the system
\[
\begin{aligned} 
\mu_{0}(u,\chi_{2},*) & = \textstyle \frac{1}{33} (-2 t_1 + 33) \geq 0; \quad \; \; 
\mu_{11}(u,\chi_{2},*) = \textstyle \frac{1}{33} (t_1 + 33) \geq 0; \\ 
\mu_{1}(u,\chi_{2},*) & = \textstyle \frac{1}{33} (-t_2 + 22) \geq 0; \qquad 
\mu_{3}(u,\chi_{2},*) = \textstyle \frac{1}{33} (2 t_2 + 22) \geq 0; \\ 
\mu_{0}(u,\chi_{3},*) & = \textstyle \frac{1}{33} (2 t_3 + 253) \geq 0; \qquad  
\mu_{11}(u,\chi_{3},*) = \textstyle \frac{1}{33} (-t_3 + 253) \geq 0; \\ 
\mu_{1}(u,\chi_{3},*) & = \textstyle \frac{1}{33} ( t_4 + 253) \geq 0; \qquad \; 
\mu_{3}(u,\chi_{3},*)  = \textstyle \frac{1}{33} (-2 t_4 + 253) \geq 0; \\ 
\mu_{0}(u,\chi_{6},*) & = \textstyle \frac{1}{33} (2 t_5 + 891) \geq 0; \qquad 
\mu_{11}(u,\chi_{6},*)  = \textstyle \frac{1}{33} (- t_5 + 891) \geq 0; \\ 
\mu_{1}(u,\chi_{6},*) & = \textstyle \frac{1}{33} (t_6 + \alpha) \geq 0; \qquad \quad \; 
\mu_{3}(u,\chi_{6},*) = \textstyle \frac{1}{33} (-2 t_6 + \alpha) \geq 0, \\ 
\end{aligned} 
\]
where the values of $\alpha$ are given in the following table

\smallskip
\centerline{\small{
\begin{tabular}{|c|c|c|c|c|c|c|c|c|c|c|c|c|}
\hline
$(\nu_{11a},\nu_{11b})$& (12, -11) &(11,-10) &(10, -9) &  (9, -8) & (8,-7) & (7,-6) & (6,-5) & (5,-4) \\ 
\hline
$\alpha$ & 1023      & 1012    & 1001 & 990   & 979  & 968  & 957  & 946  \\
\hline\hline
$(\nu_{11a},\nu_{11b})$& (4,-3) & (3,-2) &(2, -1)& (1,0) & (0,1) & (-1,2) & (-2,3) & (-3,4) \\
\hline
$\alpha$ & 935  & 924  & 913& 902 & 891 & 880 & 869 & 858 \\
\hline\hline
$(\nu_{11a},\nu_{11b})$& (-4,5) & (-5,6) & (-6,7) & (-7,8) & (-8,9) & (-9,10) & (-10,11) & (-11,12) \\
\hline
$\alpha$ & 847 & 836 & 825 & 814 & 803 & 792 & 781 & 770 \\
\hline
\end{tabular}
}}

and 
\[
\begin{aligned}
t_1 & = 40 \nu_{3a} - 50 \nu_{3b} + 10 \nu_{3c} - 10 \nu_{11a} - 10 \nu_{11b} + 
4 \nu_{3a}^{(11)} - 5 \nu_{3b}^{(11)} +  \nu_{3c}^{(11)}; \\
t_2 & = 4 \nu_{3a} - 5 \nu_{3b} +  \nu_{3c} -  \nu_{11a} -  \nu_{11b} - 4 \nu_{3a}^{(11)} + 
5 \nu_{3b}^{(11)} -  \nu_{3c}^{(11)}; \\
t_3 & = 100 \nu_{3a} + 100 \nu_{3b} + 10 \nu_{3c} + 10 \nu_{3a}^{(11)} + 
10 \nu_{3b}^{(11)} +  \nu_{3c}^{(11)};\\
t_4 & = 10 \nu_{3a} + 10 \nu_{3b} +  \nu_{3c} - 10 \nu_{3a}^{(11)} - 
10 \nu_{3b}^{(11)} -  \nu_{3c}^{(11)};\\
t_5 &= 320 \nu_{3a} - 40 \nu_{3b} - 70 \nu_{3c} - 5 \nu_{11a} - 5 \nu_{11b} + 
32 \nu_{3a}^{(11)} - 4 \nu_{3b}^{(11)} - 7 \nu_{3c}^{(11)};\\
t_6 & = 32 \nu_{3a} - 4 \nu_{3b} - 7 \nu_{3c} - 6 \nu_{11a} + 5 \nu_{11b} - 
32 \nu_{3a}^{(11)} + 4 \nu_{3b}^{(11)} + 7 \nu_{3c}^{(11)}.\\
\end{aligned}
\]

\noindent
It is not difficult to check that for any value of $\alpha$ listed in the table
this system has no integral solutions such that
each $\mu_{i}(u,\chi_{j},*)$ is a non-negative integer.

%
%

\noindent $\bullet$ Let $u$ be a unit of order $35$.
Using Proposition \ref{P:Luthar-Passi} for the 3-Brauer character $\chi_{3}$,
for which $\chi(C_{5}) = 1$ and $\chi(C_{7}) = 0$, we obtain the system
\[
\mu_{0}(u,\chi_{3},3)  = \textstyle \frac{1}{35}(24 \nu_{5} + 130) \geq 0; \quad 
\mu_{7}(u,\chi_{3},3)  = \textstyle \frac{1}{35}(-6 \nu_{5} + 125) \geq 0, \\
\]
from which $(\nu_5,\nu_7)=(15,-14)$. 

Any ordinary or Brauer $(5,7)$-constant character 
can not eliminate this pair. Indeed, 
since we have only one conjugacy class of elements of order 7,
to enumerate $(5,7)$-irreducible characters we need to look only
on character values on elements of order 5.
First, for ordinary characters and Brauer characters for $p \in \{ 11, 23 \}$
the set of differences $\chi(5a)-\chi(5b)$ is $\{ \pm 5, \pm 10, \pm 15 \}$. 

Thus, besides irreducible $(5,7)$-constant characters, all other $(5,7)$-irreducible 
characters can be parametrized by the tuples from the following set: 
\[
\begin{aligned}
\{ \; & (5,5), \; (10,-10), \; (15,-15), \; (10,-5,-5), \; (5,5,-10), \; (15,-5,-5,-5), \\
      & (5,5,5,-15), \; (15,-5,-10), \; (5,10,-15), \; (10,10,-15,-5), \\
      & (15,5,-10,-10), \; (15,15,-10,-10,-10), \; (10,10,10,-15,-15) \; \},
\end{aligned}
\]
where, for example, the tuple $(5,5,-10)$ 
mean that the character is the sum of three irreducible characters, for two of them 
$\chi(5a)-\chi(5b)=5$ and for the last one $\chi(5a)-\chi(5b)=-10$. 
Enumerating all possible combinations of characters 
for each tuple from the list above, and any other $(5,7)$-constant character would 
be a sum of some already known $(5,7)$-constant characters. We need also to repeat
the same procedure for 3-Brauer characters, where 
$\chi(5a)-\chi(5b) \in \{ \pm 5, \pm 10, 15, 20 \}$, 
and for 2-Brauer characters, where 
$\chi(5a)-\chi(5b) \in \{ \pm 5, \pm 10, -20 \}$. 
Using the GAP system \cite{GAP},
we verified that $\mu_{i}(u,\chi,p)$ are non-negative integers when $(\nu_5,\nu_7)=(15,-14)$ 
for any $(5,7)$-constant character $\chi$ obtained by the procedure described above,
so it is not possible to prove non-existence of torsion units of order 35 using this
method.

For further detalisation, we consider six cases from part (v) of Theorem \ref{T:1}.
If 
$$
\chi(u^{7}) \in \{ \; \chi(5a), \; \chi(5b), \; -\chi(5a)+2\chi(5b), \; -2\chi(5a)+3\chi(5b) \; \},
$$
then we combine the condition
$\mu_{0}(u,\chi_{2},*) = \textstyle \frac{1}{35} ( -24 t + \alpha) \geq 0$, 
where $t = 2 \nu_{5a} - 3 \nu_{5b} - 2 \nu_{7a}$
and $\alpha$ is equal to 27, 47, 67 and 87 respectively,
with the condition
$\mu_{5}(u,\chi_{2},*) = \textstyle \frac{1}{35} ( 4 t + 13) \geq 0$, 
when $\chi(u^{7}) = \chi(5a)$ and with
$\mu_{7}(u,\chi_{2},*) = \textstyle \frac{1}{35} (6 t + \beta) \geq 0$,
where $\beta$ is equal to 32, 27 and 22 respectively in the other three cases
to show that there are no solutions.

In the remaining two cases we obtain the system of inequalities
\[
\begin{aligned}
\mu_{1}(u,\chi_{2},*) & = \textstyle \frac{1}{35} (t_1 + \gamma) \geq 0; \qquad \quad 
\mu_{7}(u,\chi_{2},*)   = \textstyle \frac{1}{35} (-6t_1 + \delta) \geq 0; \\ 
\mu_{0}(u,\chi_{3},*) & = \textstyle \frac{1}{35} (24t_2 + 271) \geq 0; \quad 
\mu_{7}(u,\chi_{3},*)   = \textstyle \frac{1}{35} (-6t_2 + 256) \geq 0, \\ 
\end{aligned}
\]
where $t_1=-2 \nu_{5a} + 3 \nu_{5b} + 2 \nu_{7a}$,
$t_2=3\nu_{5a} +3\nu_{5b} +\nu_{7a}$ and
$(\gamma,\delta)=(3,17)$ for $\chi(u^{7}) = -3\chi(5a)+4\chi(5b)$
and $(\gamma,\delta)=(-2,12)$ for $\chi(u^{7}) = -4\chi(5a)+5\chi(5b)$.
Each of these cases has one solution $(\nu_{5a},\nu_{5b},\nu_{7a})=( 3, 12, -14 )$
and $( 3, 12, -14 )$ respectively.
Note that in both cases $(\nu_5,\nu_7)=(15,-14)$ as it was concluded before.

$\bullet$ It remains to prove 
that $V(\mathbb ZG)$ has no elements of orders
46, 55, 69, 77, 115, 161 and 253.
We give a detailed proof for the order 46. Other cases can be
derived similarly from the table below containing the data for
the constraints on partial augmentations $\nu_p$ and $\nu_q$
for possible orders $pq$ (including the order 46 as well)
accordingly to (\ref{E:3})--(\ref{E:5}).

If $|u|=46$, using Proposition \ref{P:Luthar-Passi} for the ordinary character 
$\chi_{23}$ with $\chi(C_{2}) = -55$ and $\chi(C_{23}) = 0$, we obtain the system
\[
\begin{aligned}
\mu_{0}(u,\chi_{23},*) & = \textstyle \frac{110}{46}(-11 \nu_{2} + 287) \geq 0; \qquad
\mu_{1}(u,\chi_{23},*)  = \textstyle \frac{55}{46}(- \nu_{2} + 576) \geq 0; \\
& \qquad \mu_{23}(u,\chi_{23},*) = \textstyle \frac{110}{46}(11 \nu_{2} + 288) \geq 0, \\
\end{aligned}
\]
which yields $\nu_2 \in \{ -22, 24 \}$. Put $t=462 \nu_{2} -22 \nu_{23}$.
Using the ordinary
character $\chi=\chi_2+\chi_5+\chi_8$ such 
that $\chi(C_{2}) = 21$, $\chi(C_{23}) = -1$
we eliminate the case $\nu_2 = -22$ from the condition
$\mu_{0}(u,\chi,*) = \textstyle \frac{1}{46}( t + 2068) \geq 0$
and the case $\nu_2 = 24$ from
$\mu_{23}(u,\chi,*) = \textstyle \frac{1}{46}(-t + 2026) \geq 0$.

The data for the remaining orders are given in the table below.

\smallskip
\centerline{\small{
\begin{tabular}{|c|c|c|c|c|c|c|c|c|c|c|c|c|c|}
\hline
$|u|$&$p$&$q$&$\xi, \; \tau$&$\xi(C_p)$&$\xi(C_q)$&$l$&$m_1$&$m_p$&$m_q$ \\
\hline
   &   &   & $\xi=(23)_{[*]}$      & -55& 0 & 1 & 31570 & -1210 & 0 \\
   &   &   & $\xi=(23)_{[*]}$      & -55& 0 & 0 & 31680 & -55   & 0 \\
46 & 2 & 23& $\xi=(23)_{[*]}$      & -55& 0 & 23& 31680 & 1210  & 0 \\
   &   &   & $\tau=(2,5,8)_{[*]}$  & 21 & -1& 0 & 2068  & 462   & -22 \\   
   &   &   & $\tau=(2,5,8)_{[*]}$  & 21 & -1& 23& 2026  & -462  & 22 \\   
\hline
   &   &   &                       &   &    & 0 & 265   & 120   & 0 \\
55 & 5 & 11& $\xi=(3)_{[*]}$       & 3 & 0  & 1 & 250   & 3     & 0 \\
   &   &   &                       &   &    & 11& 250   & -30   & 0 \\
\hline
   &   &   & $\xi=(3,5,8,15,19)_{[5]}$ & 25& 0 & 0 & 48189 & 1100 & 0 \\
69 & 3 & 23& $\xi=(3,5,8,15,19)_{[5]}$ & 25& 0 & 23& 48114 & -550 & 0 \\
   &   &   & $\tau=(3,3,6)_{[2]}$      & 12& 1 & 23& 1966  & -264 & -22 \\
\hline
77 & 7 & 11& $\xi=(3)_{[*]}$       & 1 & 0 & 0 & 259  & 60     & 0 \\
   &   &   &                       &   &   & 11& 252  & -10    & 0 \\
\hline
115& 5 & 23& $\xi=(3)_{[*]}$       & 3 & 0 & 0 & 265  & 264    & 0 \\
   &   &   &                       &   &   & 23& 250  & -66    & 0 \\
\hline
161& 7 & 23& $\xi=(2)_{[*]}$       & 2 & 0 & 0 & 35  & 264    & 0 \\
   &   &   &                       &   &   & 23& 21  & -44    & 0 \\
\hline
   &   &   &                       &   &   & 0 & 33   & 220  & 0 \\
253& 11& 23& $\xi=(2)_{[*]}$       & 1 & 0 & 11& 33   & -10  & 0 \\
   &   &   &                       &   &   & 23& 22   & -22  & 0 \\
\hline
\end{tabular}
}}
\smallskip

In the table above the notation $\xi=(i)_{[p]}$ is used for $p$-Brauer characters.

\section{Proof of Theorem \ref{T:2} }

Let $G \cong \texttt{Co}_2$. 
Then $|G|= 2^{18} \cdot 3^{6} \cdot 5^{3} \cdot 7 \cdot 11 \cdot 23$ 
and $ exp(G) = 2^{4} \cdot 3^{2} \cdot 5 \cdot 7 \cdot 11 \cdot 23$ 
(see \cite{Conway,GAP}). 
The ordinary and $p$-Brauer character tables of $G$ for
$p\in\{2,3,5,7,11,23\}$ can be found using the GAP system \cite{GAP}, 
and the same remarks as in the case of $\texttt{Co}_3$ regarding the notation for
characters and conjugacy classes applies.

The group $G$ possesses elements of orders
2, 3, 4, 5, 6, 7, 8, 9, 10, 11, 12, 14, 15, 16, 18, 20, 23, 24, 28 and 30.
As before, first we consider units of prime orders: 2, 3, 5, 7, 11 and 23.
Then we investigate products of two different primes from this list 
which are not orders of elements of $G$.

\noindent $\bullet$ Let $|u| \in \{7,11\}$. Using
Proposition \ref{P:Hertweck} we obtain that all partial
augmentations except one are zero. Thus by Proposition \ref{P:MRSW}
part (ii) of  Theorem \ref{T:2} is proved.

%
%

\noindent $\bullet$ Let $|u|=2$.
By (\ref{E:1}) and Proposition \ref{P:Hertweck} we get
$\nu_{2a}+\nu_{2b}+\nu_{2c}=1$. Put $t_1=9 \nu_{2a} - 7 \nu_{2b} +  \nu_{2c}$ and $t_2=29 \nu_{2a} + 13 \nu_{2b} - 11 \nu_{2c}$. Now using Proposition \ref{P:Luthar-Passi}
we get
\[
\begin{aligned}
\mu_{1}(u,\chi_{2},*) & = \textstyle \frac{1}{2} (t_1 + 23) \geq 0; \qquad 
\mu_{0}(u,\chi_{2},*) = \textstyle \frac{1}{2} (-t_1 + 23) \geq 0; \\ 
\mu_{0}(u,\chi_{3},*) & = \textstyle \frac{1}{2} (t_2 + 253) \geq 0; \quad 
\mu_{1}(u,\chi_{3},*)  = \textstyle \frac{1}{2} (-t_2 + 253) \geq 0, \\ 
\end{aligned}
\]
which has only 48 solutions listed in the Appendix B
such that all $\mu_{i}(u,\chi_{j},*)$ are non-negative integers.

%
%

\noindent $\bullet$ Let $u$ be a unit of order $3$. By (\ref{E:1})
and Proposition \ref{P:Hertweck} we get $\nu_{3a}+\nu_{3b}=1$. Put $t_1=4 \nu_{3a} - 5 \nu_{3b}$. Using
Proposition \ref{P:Luthar-Passi} we obtain the  following system
$$
\mu_{0}(u,\chi_{2},*) = \textstyle \frac{1}{3} (-2t_1 + 23) \geq 0; \qquad 
\mu_{1}(u,\chi_{2},*) = \textstyle \frac{1}{3} (t_1 + 23) \geq 0, \\ 
$$
which has only four solutions listed in part (iv) of Theorem \ref{T:2}
such that all $\mu_{i}(u,\chi_{2},*)$ are non-negative integers.

%
%

\noindent $\bullet$ Let $u$ be a unit of order $5$.
By (\ref{E:1}) and Proposition \ref{P:Hertweck} we get
$\nu_{5a}+\nu_{5b}=1$. If put $t_1=2 \nu_{5a} - 3 \nu_{5b}$, then using Proposition \ref{P:Luthar-Passi}
we obtain 
$$
\mu_{0}(u,\chi_{2},*) = \textstyle \frac{1}{5} (-4t_1 + 23) \geq 0; \qquad 
\mu_{1}(u,\chi_{2},*) = \textstyle \frac{1}{5} (t_1 + 23) \geq 0, \\ 
$$
which has only six solutions listed in part (v) of Theorem \ref{T:2}
such that all $\mu_{i}(u,\chi_{2},*)$ are non-negative integers.

%
%

\noindent $\bullet$ Let $u$ be a unit of order $23$. By
(\ref{E:1}) and Proposition \ref{P:Hertweck} we get
$\nu_{23a}+\nu_{23b}=1$. Put $t_1=12 \nu_{23a} - 11 \nu_{23b}$.
Using Proposition \ref{P:Luthar-Passi} we obtain the system
\[
\begin{aligned}
\mu_{1}(u,\chi_{4},2) & = \textstyle \frac{1}{23} (-t_1 + 748)
\geq 0; \qquad \mu_{1}(u,\chi_{9},3) = \textstyle \frac{1}{23} (t_1+ 9372) \geq 0; \\ 
& \qquad \mu_{5}(u,\chi_{4},2) = \textstyle \frac{1}{23} (11 \nu_{23a} - 12 \nu_{23b} + 748) \geq 0, \\ 
\end{aligned}
\]
which has only 66 solutions listed in part (vi) of Theorem \ref{T:2}
such that all $\mu_{i}(u,\chi_{j})$ are non-negative integers.

Now we will deal with torsion units of orders which do not appear in $G$.

%
%

\noindent $\bullet$ Let $|u|=21$.
We need to consider four cases defined by part (iv)
of Theorem \ref{T:2}.

Case 1.
$\chi(u^{7}) = \chi(3a)$. Put $t_1 = 4 \nu_{3a} - 5 \nu_{3b} - 2 \nu_{7a}$. 
Then we obtain the system
\[
\mu_{3}(u,\chi_{2},*) = \textstyle \frac{1}{21} (2t_1 + 13) \geq 0; \qquad 
\mu_{0}(u,\chi_{2},*)  = \textstyle \frac{1}{21} (-12t_1 + 27) \geq 0, \\ 
\]
which has no solution. 

Case 2.
$\chi(u^{7}) = \chi(3b)$. Put $t_1=-4 \nu_{3a} + 5 \nu_{3b} + 2 \nu_{7a}$, $t_2=10 \nu_{3a} + 10 \nu_{3b} +  \nu_{7a}$. Then we obtain the following incompatible system of inequalities:
\[
\begin{aligned}
\mu_{0}(u,\chi_{2},*) & = \textstyle \frac{1}{21} (12t_1+ 45) \geq 0; \quad \;
\mu_{1}(u,\chi_{2},*)  = \textstyle \frac{1}{21} (t_1+ 16) \geq 0; \\ 
\mu_{1}(u,\chi_{3},*) & = \textstyle \frac{1}{21} (t_2+ 242) \geq 0; \qquad
\mu_{7}(u,\chi_{3},*)  = \textstyle \frac{1}{21} (-6t_2 + 249) \geq 0; \\ 
\mu_{7}(u,\chi_{2},*) & = \textstyle \frac{1}{21} (-6 t_1 + 30) \geq 0; \quad 
\mu_{0}(u,\chi_{3},*)  = \textstyle \frac{1}{21} (12 t_2 + 279) \geq 0; \\ 
&\mu_{1}(u,\chi_{5},*)  = \textstyle \frac{1}{21} (-11 \nu_{3a} + 16 \nu_{3b} + 1755) \geq 0. \\ 
\end{aligned}
\]

Case 3.
$\chi(u^{7}) = -\chi(3a)+2\chi(3b)$. Put $t_1=4 \nu_{3a} -5 \nu_{3b} - 2 \nu_{7a}$. 
The system
\[
\begin{aligned}
\mu_{7}(u,\chi_{2},*) & = \textstyle \frac{1}{21} (6t_1 + 21) \geq 0; \qquad 
\mu_{0}(u,\chi_{2},*)  = \textstyle \frac{1}{21} (-12t_1 + 63) \geq 0; \\ 
&\qquad \mu_{1}(u,\chi_{2},*)  = \textstyle \frac{1}{21} (-t_1 + 7) \geq 0, \\ 
\end{aligned}
\]
has no integral solutions.

Case 4.
$\chi(u^{7}) = -2\chi(3a)+3\chi(3b)$. Put \quad $t_1=-4 \nu_{3a} + 5 \nu_{3b} + 2 \nu_{7a}$\quad  and $t_2=10 \nu_{3a} + 10 \nu_{3b} +  \nu_{7a}$. 
Then the following system is incompatible:
\[
\begin{aligned}
\mu_{1}(u,\chi_{2},*) & = \textstyle \frac{1}{21} (t_1 - 2) \geq 0; \qquad \quad
\mu_{7}(u,\chi_{2},*) = \textstyle \frac{1}{21} (-6t_1 + 12) \geq 0; \\ 
\mu_{0}(u,\chi_{3},*) & = \textstyle \frac{1}{21} (12t_2 + 279) \geq 0; \quad 
\mu_{1}(u,\chi_{3},*)  = \textstyle \frac{1}{21} (t_2 + 242) \geq 0; \\ 
& \qquad \mu_{7}(u,\chi_{3},*) = \textstyle \frac{1}{21} (-6 t_2 + 249) \geq 0. \\ 
\end{aligned}
\]

%
%

\noindent $\bullet$ Let $u$ be a unit of order $22$. We will use the same
approach as in the case of units of order 33 in the proof of Theorem \ref{T:1}.
Taking into account part (ii) of the Theorem \ref{T:2} 
with respect to torsion units of order 11,
we obtain the system of inequalities:
\[
\begin{aligned} 
\mu_{0}(u,\chi_{2},*) & = \textstyle \frac{1}{22} (-t_1 + 33) \geq 0; \qquad  
\mu_{11}(u,\chi_{2},*)  = \textstyle \frac{1}{22} (t_1 + 33) \geq 0; \\ 
\mu_{1}(u,\chi_{2},*) & = \textstyle \frac{1}{22} (-t_2 + 22) \geq 0; \qquad 
\mu_{2}(u,\chi_{2},*)  = \textstyle \frac{1}{22} (t_2 + 22) \geq 0; \\ 
\mu_{0}(u,\chi_{3},*) & = \textstyle \frac{1}{22} (t_3 + 253) \geq 0; \qquad 
\mu_{11}(u,\chi_{3},*)  = \textstyle \frac{1}{22} (-t_3 + 253) \geq 0; \\ 
\mu_{1}(u,\chi_{3},*) & = \textstyle \frac{1}{22} (t_4 + 253) \geq 0; \qquad 
\mu_{2}(u,\chi_{3},*)  = \textstyle \frac{1}{22} (-t_4 + 253) \geq 0; \\ 
\mu_{0}(u,\chi_{4},*) & = \textstyle \frac{1}{22} (t_5 + 275) \geq 0; \qquad 
\mu_{11}(u,\chi_{4},*)  = \textstyle \frac{1}{22} (-t_5 + 275) \geq 0, \\ 
\end{aligned} 
\]
where
\[
\begin{aligned} 
t_1 &= 90 \nu_{2a} - 70 \nu_{2b} + 10 \nu_{2c} - 10 \nu_{11a} + 9 \nu_{2a}^{(11)} - 7 \nu_{2b}^{(11)} +  \nu_{2c}^{(11)}; \\
t_2 &= 9 \nu_{2a} - 7 \nu_{2b} +  \nu_{2c} -  \nu_{11a} - 9 \nu_{2a}^{(11)} + 7 \nu_{2b}^{(11)} -  \nu_{2c}^{(11)}; \\
t_3 &= 290 \nu_{2a} + 130 \nu_{2b} - 110 \nu_{2c} + 29 \nu_{2a}^{(11)} + 13 \nu_{2b}^{(11)} - 11 \nu_{2c}^{(11)} ;\\
t_4 &= 29 \nu_{2a} + 13 \nu_{2b} - 11 \nu_{2c} - 29 \nu_{2a}^{(11)} - 13 \nu_{2b}^{(11)} + 11 \nu_{2c}^{(11)} ; \\
t_5 &= 510 \nu_{2a} + 350 \nu_{2b} + 110 \nu_{2c} + 51 \nu_{2a}^{(11)} + 35 \nu_{2b}^{(11)} + 11 \nu_{2c}^{(11)},\\
\end{aligned} 
\]
which has no integral solutions such that each $\mu_{i}(u,\chi_{j},*)$
is a non-negative integer.

%
%

\noindent $\bullet$ Let $u$ be a unit of order $33$.
Using Proposition \ref{P:Luthar-Passi} for the ordinary character $\chi_3$, for which
$\chi(C_{3}) = 10$ and $\chi(C_{11}) = 0$, we obtain the system
$$
\mu_{0}(u,\chi_{3},*)  = \textstyle \frac{1}{33}(200 \nu_{3} + 273) \geq 0; \qquad
\mu_{11}(u,\chi_{3},*) = \textstyle \frac{1}{33}(-100 \nu_{3} + 243) \geq 0, \\
$$
which has no solutions such that all $\mu_{i}(u,\chi_{3},*)$ are non-negative integers.

%
%

\noindent $\bullet$ Let $u$ be a unit of order $35$.
For this order we have situation similar to the 3rd Conway group.
It is possible to show that the best restriction that we can get
is the system obtained by Proposition \ref{P:Luthar-Passi} for the 3-Brauer
character $\chi_{7}$, for which $\chi(C_{5}) = 4$ and $\chi(C_{7}) = 0$:
$$
\mu_{0}(u,\chi_{7},3) = \textstyle \frac{1}{35}(96 \nu_{5} + 2270) \geq 0; \qquad
\mu_{7}(u,\chi_{7},3) = \textstyle \frac{1}{35}(-24 \nu_{5} + 2250) \geq 0, \\
$$
from which $\nu_5 \in \{ -20, 15, 50, 85 \}$. Using further analysis
we are able to eliminate three of these opportunities, but not all.
Now considering six cases defined by part (v) of Theorem \ref{T:1},
we obtain exactly the same systems of inequalities that either has
no solutions or lead to solutions listed in part (vii) of Theorem \ref{T:2}.

$\bullet$ To complete the proof of part (i), it remains
to show that there are no elements of orders 46, 55, 69, 77, 115, 161 and 253
in $V(\mathbb Z G)$.
As in the proof of Theorem \ref{T:1}, below we give the table containing
the data describing the constraints on partial augmentations
$\nu_p$ and $\nu_q$ accordingly to (\ref{E:3})--(\ref{E:5}) for
all these orders. From this table
part (i) of Theorem \ref{T:2} is derived in the same way as in the 
proof of Theorem \ref{T:1}.

\smallskip
\centerline{\small{
\begin{tabular}{|c|c|c|c|c|c|c|c|c|c|c|c|c|c|}
\hline
$|u|$&$p$&$q$&$\xi, \; \tau$&$\xi(C_p)$&$\xi(C_q)$&$l$&$m_1$&$m_p$&$m_q$ \\
\hline
   &   &   & $\xi=(2,3,5)_{[*]}$   & -1 & 0 & 1 & 2046 & -22   & 0 \\
   &   &   & $\xi=(2,3,5)_{[*]}$   & -1 & 0 & 2 & 2046 & 1     & 0 \\
46 & 2 & 23& $\xi=(2,3,5)_{[*]}$   & -1 & 0 & 23& 2048 & 22    & 0 \\
   &   &   & $\tau=(2,4,5)_{[*]}$  & 21 & -1& 0 & 2068 & 462   & -22 \\   
   &   &   & $\tau=(2,4,5)_{[*]}$  & 21 & -1& 23& 2026 & -462  & 22 \\   
\hline
   &   &   &                       &   &    & 0 & 265  & 120   & 0 \\
55 & 5 & 11& $\xi=(3)_{[*]}$       & 3 & 0  & 5 & 265  & -12   & 0 \\
   &   &   &                       &   &    & 11& 250  & -30   & 0 \\
\hline
69 & 2 & 23& $\xi=(3)_{[*]}$       & 10& 0 & 0 & 273  & 440    & 0 \\
   &   &   &                       &   &   & 23& 243  & -220   & 0 \\
\hline
77 & 7 & 11& $\xi=(4)_{[*]}$       & 3 & 0 & 0 & 287  & 120    & 0 \\
   &   &   &                       &   &   & 11& 273  & -20    & 0 \\
\hline
115& 5 & 23& $\xi=(3)_{[*]}$       & 3 & 0 & 0 & 265  & 264    & 0 \\
   &   &   &                       &   &   & 23& 250  & -66    & 0 \\
\hline
161& 7 & 23& $\xi=(2)_{[*]}$       & 2 & 0 & 0 & 35  & 264    & 0 \\
   &   &   &                       &   &   & 23& 21  & -44    & 0 \\
\hline
   &   &   &                       &   &   & 0 & 33   & 220  & 0 \\
253& 11& 23& $\xi=(2)_{[*]}$       & 1 & 0 & 11& 33   & -10  & 0 \\
   &   &   &                       &   &   & 23& 22   & -22  & 0 \\
\hline
\end{tabular}
}}

\section{Proof of Theorem \ref{T:3} }

Let $G \cong \texttt{Co}_1$. 
Then (see \cite{Conway,GAP}) we know that
$|G|= 2^{21} \cdot 3^{9} \cdot 5^{4} \cdot 7^{2} \cdot 11 \cdot 13 \cdot 23$  
and $ exp(G) =  2^{4} \cdot 3^{2} \cdot 5 \cdot 7 \cdot 11 \cdot 13 \cdot 23$.  
The ordinary and $p$-Brauer character tables of $G$ for
$p\in\{7,11,23\}$ can be found using the GAP system \cite{GAP},
and we use the same approach to the notation for
characters and conjugacy classes. 

\noindent $\bullet$ Let $|u| \in \{11,12\}$. Using
Proposition \ref{P:Hertweck} we obtain that all partial
augmentations except one are zero. Thus by Proposition \ref{P:MRSW}
part (ii) of  Theorem \ref{T:3} is proved.

\noindent $\bullet$ Let $|u|=7$. By (\ref{E:1}) and Proposition \ref{P:Hertweck} we get
$\nu_{7a}+\nu_{7b}=1$. Put $t=10\nu_{7a} + 3 \nu_{7b}$. 
Then using Proposition \ref{P:Luthar-Passi} we obtain the system of inequalities
$$ 
\mu_{0}(u,\chi_{2},*) = \textstyle \frac{1}{7} (6 t + 276) \geq 0; \quad 
\mu_{1}(u,\chi_{2},*) = \textstyle \frac{1}{7} (-t + 276) \geq 0, \\ 
$$
which yields only 47 solutions listed in part (iii) of Theorem \ref{T:3}.

\noindent $\bullet$ Let $|u|=23$. Similarly to the above,
$\nu_{23a}+\nu_{23b}=1$ and we get the system 
\[
\begin{aligned} 
\mu_{1}(u,\chi_{17},*) & = \textstyle \frac{1}{23} (12 \nu_{23a} - 11 \nu_{23b} + 
673750) \geq 0; \\ 
\mu_{5}(u,\chi_{17},*) & = \textstyle \frac{1}{23} (-11 \nu_{23a} + 12 \nu_{23b} + 
673750) \geq 0, \\ 
\end{aligned} 
\]
which has only 58588 solutions listed in part (iv) of Theorem \ref{T:3}.

Now we consider torsion units of orders which do not appear in $G$.
First we will show that $V(\mathbb ZG)$ has no units of order
$46$, $69$, $77$, $91$, $115$, $143$, $161$, $253$ and $299$.

As in the proof of Theorem \ref{T:1}, below we give the table containing
the data describing the constraints on partial augmentations
$\nu_p$ and $\nu_q$ accordingly to (\ref{E:3})--(\ref{E:5}) for
orders $91$, $143$, $161$, $253$ and $299$.

\smallskip
\centerline{\small{
\begin{tabular}{|c|c|c|c|c|c|c|c|c|c|c|c|c|c|}
\hline
$|u|$&$p$&$q$&$\xi, \; \tau$&$\xi(C_p)$&$\xi(C_q)$&$l$&$m_1$&$m_p$&$m_q$ \\
\hline
91 & 7 & 13& $\xi=(3)_{[*]}$       & 5 & 0  & 0 & 329  & 360   & 0 \\
   &   &   &                       &   &    & 13& 294  & -60   & 0 \\
\hline
143& 2 & 0 & $\xi=(3)_{[*]}$       & 2 & 0  & 0 & 319  & 240   & 0 \\
   &   &   &                       &   &    & 13& 297  & -24   & 0 \\
\hline
161& 7 & 23& $\xi=(3)_{[*]}$       & 5 & 0  & 0 & 329  & 660   & 0 \\
   &   &   &                       &   &    & 23& 294  & -110  & 0 \\
\hline
   &   &   &                       &   &   & 0 & 319   & 440   & 0 \\
253& 11& 23& $\xi=(3)_{[*]}$       & 2 & 0 & 11& 319   & -20   & 0 \\
   &   &   &                       &   &   & 23& 297   & -44   & 0 \\
\hline
299& 13& 23& $\xi=(2)_{[*]}$       & 3 & 0  & 0 & 312  & 792   & 0 \\
   &   &   &                       &   &    & 23& 273  & -66   & 0 \\
\hline
\end{tabular}
}}
\bigskip

\noindent $\bullet$ Let $|u|=77$. Similarly to the case of order 33 in Theorem \ref{T:1}
we have that
\[
\begin{aligned} 
\mu_{11}(u,\chi_{2},*) & = \textstyle \frac{1}{77} (-t_1 - 30 \nu_{7b} - 10 \nu_{11a} + 10 \nu_{11a}^{(7)} - 3 \nu_{7b}^{(11)} + 276) \geq 0; \\ 
\mu_{0}(u,\chi_{3},*) & = \textstyle \frac{1}{77} (3 t_1 + 300 \nu_{7b} + 120 \nu_{11a} + 20 \nu_{11a}^{(7)} + 30 \nu_{7b}^{(11)} + 299) \geq 0; \\ 
\mu_{0}(u,\chi_{4},*) & = \textstyle \frac{1}{77} (6 t_2 + 1771) \geq 0; \quad \; 
\mu_{11}(u,\chi_{4},*) = \textstyle \frac{1}{77} (-t_2 + 1771) \geq 0; \\ 
\mu_{1}(u,\chi_{4},*) & = \textstyle \frac{1}{77} (t_3 + 1771) \geq 0; \qquad 
\mu_{7}(u,\chi_{4},*)   = \textstyle \frac{1}{77} (-6 t_3 + 1771) \geq 0; \\ 
\mu_{0}(u,\chi_{7},*) & = \textstyle \frac{1}{77} (2 t_4 + 27300) \geq 0; \quad 
\mu_{0}(u,\chi_{15},13) = \textstyle \frac{1}{77} (-t_4 + 474145) \geq 0, \\ 
\end{aligned} 
\]
where $t_1=100 \nu_{7a} + 10 \nu_{7a}^{(11)}$, $t_2= 140 \nu_{7a} + 14 \nu_{7a}^{(11)}$,
$t_3= 14 \nu_{7a} - 14 \nu_{7a}^{(11)}$ and 
$t_4=420 \nu_{7a} - 60 \nu_{11a} - 10 \nu_{11a}^{(7)} + 42 \nu_{7a}^{(11)}$. 
This system has no integral solutions such that all $\mu_{i}(u,\chi_{j},p)$
are non-negative integers.

\noindent $\bullet$ Let $|u| \in \{46, 69, 115\}$. Using the same method as in the case
of order 33 in Theorem \ref{T:1}, we constructed systems of constraints 
and verified with the constraint solver ECLiPSe \cite{ECLiPSe} that they 
have no solutions, using the lower and upper bounds on partial augmentations
given by Proposition \ref{P:Hales-Luthar-Passi}. Remarkably, we immediately 
check that there are no units of order $69$, while the enumeration of 
all possible partial augmentations for $|u| \in \{3,23\}$
requires $15239 \cdot 58588 = 892822532$ cases.

\noindent $\bullet$ Let $|u| \in \{55, 65 \}$. Clearly, $|u^{11}|=|u^{13}|=5$. 
Using the LAGUNA package \cite{LAGUNA} together
with ECLiPSe \cite{ECLiPSe}, we produce 1041 possible tuples
of partial augmentations for units of order 5. From 
corresponding systems we computed solutions
listed in parts (v) and (vi). 
Note that $p$-Brauer character tables for $G$ are 
are not known for $p \in \{2,3,5\}$
(see \verb+http://www.math.rwth-aachen.de/~MOC/work.html+), and
hopefully, further progress could be made 
when $p$-Brauer character tables for missing values of $p$ will became available.

\section*{Acknowledgments}
The research was supported by OTKA No.K68383.

\bibliographystyle{plain}
\bibliography{conway}

\newpage

\centerline{\bf Appendix A.}

\vspace{5pt}

Partial augmentations $(\nu_{3a},\nu_{3b},\nu_{3c})$ for units of
order 3 in $\mathbb Z\texttt{Co}_3$:
$$
{\tiny
\begin{array}{llllll}
( -9, -3, 13 ), & ( -9, -2, 12 ), & ( -9, -1, 11 ), & ( -8, -5, 14 ), & ( -8, -4, 13 ), & ( -8, -3, 12 ), \\
( -8, -2, 11 ), & ( -8, -1, 10 ), & ( -8, 0, 9 ),   & ( -7, -5, 13 ), & ( -7, -4, 12 ), & ( -7, -3, 11 ), \\
( -7, -2, 10 ), & ( -7, -1, 9 ), & ( -7, 0, 8 ), & ( -6, -4, 11 ), & ( -6, -3, 10 ), & ( -6, -2, 9 ), \\
( -6, -1, 8 ),  & ( -6, 0, 7 ), & ( -6, 1, 6 ), & ( -5, -4, 10 ), & ( -5, -3, 9 ),  & ( -5, -2, 8 ), \\
( -5, -1, 7 ), & ( -5, 0, 6 ),   & ( -5, 1, 5 ),  & ( -4, -3, 8 ), & ( -4, -2, 7 ),  & ( -4, -1, 6 ), \\
( -4, 0, 5 ), & ( -4, 1, 4 ), & ( -4, 2, 3 ),   & ( -3, -3, 7 ),  & ( -3, -2, 6 ), & ( -3, -1, 5 ), \\
( -3, 0, 4 ),   & ( -3, 1, 3 ), & ( -3, 2, 2 ), & ( -2, -2, 5 ), & ( -2, -1, 4 ), & ( -2, 0, 3 ), \\
( -2, 1, 2 ), & ( -2, 2, 1 ), & ( -2, 3, 0 ), & ( -1, -2, 4 ), & ( -1, -1, 3 ), & ( -1, 0, 2 ), \\
( -1, 1, 1 ), & ( -1, 2, 0 ), & ( -1, 3, -1 ), & ( 0, -1, 2 ), & ( 0, 0, 1 ), & ( 0, 1, 0 ), \\
( 0, 2, -1 ), & ( 0, 3, -2 ), & ( 0, 4, -3 ), & ( 1, -1, 1 ), & ( 1, 0, 0 ), & ( 1, 1, -1 ), \\
( 1, 2, -2 ), & ( 1, 3, -3 ), & ( 1, 4, -4 ), & ( 2, 0, -1 ), & ( 2, 1, -2 ), & ( 2, 2, -3 ), \\
( 2, 3, -4 ), & ( 2, 4, -5 ), & ( 2, 5, -6 ), & ( 3, 0, -2 ), & ( 3, 1, -3 ), & ( 3, 2, -4 ), \\
( 3, 3, -5 ), & ( 3, 4, -6 ), & ( 3, 5, -7 ), & ( 4, 1, -4 ), & ( 4, 2, -5 ), & ( 4, 3, -6 ), \\
( 4, 4, -7 ), & ( 4, 5, -8 ), & ( 4, 6, -9 ), & ( 5, 1, -5 ), & ( 5, 2, -6 ), & ( 5, 3, -7 ), \\
( 5, 4, -8 ), & ( 5, 5, -9 ), & ( 5, 6, -10 ), & ( 6, 2, -7 ), & ( 6, 3, -8 ), & ( 6, 4, -9 ), \\
( 6, 5, -10 ), & ( 6, 6, -11 ), & ( 6, 7, -12 ), & ( 7, 2, -8 ), & ( 7, 3, -9 ), & ( 7, 4, -10 ), \\
( 7, 5, -11 ), & ( 7, 6, -12 ), & ( 7, 7, -13 ), & ( 8, 3, -10 ), & ( 8, 4, -11 ), & ( 8, 5, -12 ), \\
( 8, 6, -13 ), & ( 8, 7, -14 ), & ( 8, 8, -15 ), & ( 9, 3, -11 ), & ( 9, 4, -12 ), & ( 9, 5, -13 ), \\
( 9, 6, -14 ), & ( 9, 7, -15 ), & ( 9, 8, -16 ), & ( 10, 4, -13 ), & ( 10, 5, -14 ), & ( 10, 6, -15 ), \\
( 10, 7, -16 ), & ( 10, 8, -17 ), & ( 10, 9, -18 ), & ( 11, 4, -14 ), & ( 11, 5, -15 ), & ( 11, 6, -16 ), \\
( 11, 7, -17 ), & ( 11, 8, -18 ), & ( 11, 9, -19 ), & ( 12, 5, -16 ), & ( 12, 6, -17 ), & ( 12, 7, -18 ), \\
( 12, 8, -19 ), & ( 12, 9, -20 ), & ( 12, 10, -21 ), & ( 13, 5, -17 ), & ( 13, 6, -18 ), & ( 13, 7, -19 ), \\
( 13, 8, -20 ), & ( 13, 9, -21 ), & ( 13, 10, -22 ), & ( 14, 6, -19 ), & ( 14, 7, -20 ), & ( 14, 8, -21 ), \\
( 14, 9, -22 ), & ( 14, 10, -23 ), & ( 14, 11, -24 ), & ( 15, 6, -20 ), & ( 15, 7, -21 ), & ( 15, 8, -22 ), \\
( 15, 9, -23 ), & ( 15, 10, -24 ), & ( 15, 11, -25 ), & ( 16, 7, -22 ), & ( 16, 8, -23 ), & ( 16, 9, -24 ), \\
( 16, 10, -25 ), & ( 16, 11, -26 ), & ( 16, 12, -27 ), & ( 17, 7, -23 ), & ( 17, 8, -24 ).\\
\end{array}
}
$$

\vspace{5pt}

\centerline{\bf Appendix B.}

\vspace{5pt}

Partial augmentations $(\nu_{2a},\nu_{2b},\nu_{2c})$ for units of
order 2 in $\mathbb Z\texttt{Co}_2$:
$$
{\tiny
\begin{array}{lllllll}
( -4, -3, 8 ), & ( -4, -2, 7 ), & ( -4, -1, 6 ), & ( -3, -5, 9 ), & ( -3, -4, 8 ), & ( -3, -3, 7 ), & ( -3, -2, 6 ), \\
( -3, -1, 5 ), & ( -3, 0, 4 ), & ( -2, -4, 7 ), & ( -2, -3, 6 ), & ( -2, -2, 5 ), & ( -2, -1, 4 ), & ( -2, 0, 3 ), \\
( -2, 1, 2 ), & ( -1, -3, 5 ), & ( -1, -2, 4 ), & ( -1, -1, 3 ), & ( -1, 0, 2 ), & ( -1, 1, 1 ), & ( -1, 2, 0 ), \\
( 0, -2, 3 ), & ( 0, -1, 2 ), & ( 0, 0, 1 ), & ( 0, 1, 0 ), & ( 0, 2, -1 ), & ( 0, 3, -2 ), & ( 1, -1, 1 ), \\
( 1, 0, 0 ), & ( 1, 1, -1 ), & ( 1, 2, -2 ), & ( 1, 3, -3 ), & ( 1, 4, -4 ), & ( 2, 0, -1 ), & ( 2, 1, -2 ), \\
( 2, 2, -3 ), & ( 2, 3, -4 ), & ( 2, 4, -5 ), & ( 2, 5, -6 ), & ( 3, 1, -3 ), & ( 3, 2, -4 ), & ( 3, 3, -5 ), \\
( 3, 4, -6 ), & ( 3, 5, -7 ), & ( 3, 6, -8 ), & ( 4, 2, -5 ), & ( 4, 3, -6 ), & ( 4, 4, -7 ). \\
\end{array}
}
$$

\end{document}